\documentclass[12pt]{article}
\usepackage{a4wide,amsmath,amssymb, amsthm, amsfonts}
\usepackage[utf8]{inputenc}
\usepackage[english]{babel}
\usepackage{xcolor}
\usepackage{hyperref}
\usepackage{enumitem}
\usepackage{graphicx}
\usepackage{hyperref}
\usepackage{mathtools}
\usepackage{bbm}
\usepackage{csquotes}
\usepackage{comment}

\newcommand{\blue}[1]{\textcolor{blue}{#1}}

  \newtheorem{theo}{Theorem}[section]

  \newtheorem{proposition}[theo] {Proposition}
  
  \newtheorem{auxiliary}[theo] {Auxiliary result}

  \newtheorem{remark}[theo] {Remark}
  
  \newtheorem{assumption}[theo] {Assumption}

\def\R{\mathbb{R}}

\def\cQ{\mathcal{Q}}
\def\cU{\mathcal{U}}

\def\cC{\mathcal{C}}

\def\cL{\mathcal{L}}
\def\F{\mathbf{F}}
\def\G{\mathbf{G}}
\def\Ccoe{C_{\text{coe}}}
\def\I{[0,\infty)}
\def\T{[0,T)}

\def\util{\tilde{u}}
\def\Mtil{\widetilde{M}}
\def\Ntil{\widetilde{N}}
\def\Qtil{\widetilde{Q}}
\def\errpar{e}
\def\errst{r}

\def\lan{\langle}
\def\ran{\rangle}
\def\dupair#1{\lan #1 \ran}

\colorlet{blue}{black}

\begin{document}

\title{A model reference adaptive system approach for nonlinear online parameter identification}

\date{}

\author{Barbara Kaltenbacher, Tram Thi Ngoc Nguyen}

\author{Barbara Kaltenbacher \thanks{Department of Mathematics, Alpen-Adria-Universit\"at Klagenfurt, 
Universit\"atsstra\ss e 65-67, 9020 Klagenfurt, Austria ({\tt barbara.kaltenbacher@aau.at})} \and Tram Thi Ngoc Nguyen\thanks{Institute of Mathematics and Scientific Computing, University of Graz, Heinrichstraße 36, A-8010 Graz, Austria ({\tt tram.nguyen@uni-graz.at}). Corresponding author.}  }

\maketitle

{\bf Abstract.}
Dynamical systems, for instance in model predictive control, often contain unknown parameters, which must be determined during system operation. Online, or on-the-fly, parameter identification methods are therefore necessary. The challenge of online methods is that one must continuously estimate parameters as experimental data becomes available. The existing techniques in the context of time-dependent partial differential equations exclude the case where the system depends nonlinearly on the parameters. Based on a model reference adaptive system approach, we present an online parameter identification method for nonlinear  infinite-dimensional evolutionary systems.\\


\section{Introduction}\label{sec::intro}
Evolutionary, spatially dependent processes in science and engineering  are usually modeled by time dependent partial differential equations (PDEs). These models often contain finite or infinite dimensional parameters, such as spatially varying drift or diffusion coefficients, as well as source terms, whose values are unknown and have to be recovered on the basis of indirect observations of the system. Whenever parameter identification needs to take place during the operation of the considered system -- as is the case in e.g. model predictive control -- so-called online methods have to be employed.

Model reference adaptive systems (MRAS) do so by setting up a dynamic update law for both the parameter and the state, where the evolutionary system for the state is a modification of the original model, while the parameter evolution is driven by the observation mismatch. Additionally, stabilizing terms are introduced.
Online parameter identification has been extensively studied in the finite dimensional setting, e.g. \cite{Ioannou,Narendra05,Sastry}; however, when it come to infinite dimensional models as arising in the context of PDEs, less is known so far. 
We refer to the literature review in \cite{Baumeister-etal,Kugler1,Kugler} for MRAS-based approaches in the PDE context. In the field of inverse problems \cite{BarbaraNeubauerScherzer,Kirsch}, the idea of treating the unknown parameter and the state as an extended variable can be founded in \cite{HaAs01, aao16, Kaltenbacher, Nguyen} Following up on \cite{BoigerKaltenbacher}, which in turn was strongly inspired by \cite{Baumeister-etal} and \cite{Kugler}, we here extend the scope to problems that are nonlinear not only with respect to the state, but also with respect to the parameter -- a situation that could not be tackled by the MRAS approaches investigated so far.
To this end, we first focus on the situation of complete state observations. This may be extended to the case of partial observations, following the principles of \cite{BoigerKaltenbacher, Kugler}, however with additional  difficulties and obstacles arising due to parameter nonlinearity. As such, investigations on partial observations are postponed to future work.
In this context, we wish to point to recent progress on ensemble Kalman filters for inverse problems, see, e.g., \cite{Albers:2019,DihlmannHaasdonk2015,SchillingsStuart2017}, which are clearly promising in the online identification setting. However, note that this approach, being based on statistical considerations, is quite different from the one we are following here.

The paper is structured as follows.
In Section \ref{sec::intro}, we first introduce the model class under consideration. Then, the adaptive system that we propose for the online parameter identification task is presented.
Section \ref{sec:Error} provides an error analysis and convergence results for both exact and noisy data. Moreover, well-posedness of the adaptive system and therewith well-definedness of the method is established there.
Finally, in Section \ref{sec:ex} we discuss two examples of coefficient identification problems in parabolic PDEs, where the convergence conditions from Section \ref{sec:Error} can be verified. The Appendix supplements some auxiliary results.

\subsection{\blue{The underlying state system}}\label{sec::underlying_system}
We consider the problem of identifying the stationary parameter $q^*\in Q$ in the evolution equation
\begin{align}
& D_tu^*(t) + f(q^*,u^*(t)) = g(t) \qquad t>0 \label{model-1}\\
& u^*(0) = u_0 \label{model-2}
\end{align}
from full observations $z$ of the state $u^*$ over time
\begin{align}
z(t)=Gu^*(t)=u^*(t) \qquad t>0. \label{model-3}
\end{align}
\blue{Here, the solution $u^*$, also called the exact state, is a function on the infinite time interval $[0,\infty)$ and some bounded domain $\Omega$ with smooth boundary. The model operator $f$ is nonlinear in both $u^*$ and $q^*$. Whenever \eqref{model-1} is a PDE, $f$ will be a differential operator with respect to some spatial variables acting on the state $u^*$, with possibly space and time dependent coefficients (cf. Remark \ref{rem:non-autonomous} below) equipped with (e.g. Dirichlet or Neumann) boundary conditions; $q^*$ then typical plays the role of some space dependent coefficient or source term.}

We introduce the triples
\begin{align*}
& Q\hookrightarrow U \hookrightarrow Z,\quad U\hookrightarrow H\hookrightarrow U^* , \quad Q\hookrightarrow H\hookrightarrow Q^*
\end{align*}
of function spaces, where the last two are Gelfand triples. \blue{Here, $\hookrightarrow$ refers to continuous embeddings between Banach spaces \cite{Leoni:2009}. In this setting, $Q$ (parameter space), $U$ (state space), $Z$ (observation space) are reflexive Banach spaces, and $H$ is a Hilbert space}.
Furthermore, \blue{by a usual overload of} notation, we denote  by $f$ also
the Nemytskii operator 
\begin{align*}
&f:L^2(\I;H)\times L^2(\I;U)\to L^2(\I;U^*)\\ 
& [f(q,u)](t)=f(q(t),u(t))
\end{align*}
induced by $f:H\times U \to U^*$.
In addition to the assumption 
\begin{align*}
& g\in L^2(\I;U^*)\\
& u_0 \in H,
\end{align*}
on the source term and initial data, we suppose that
\begin{align*}
& G=\text{Id}:U\to Z(\supseteq U),
\end{align*}
\blue{meaning that the observation operator $G$ is a continuous embedding from $U$ to $Z$, resulting  in $z$ being a full measurement in time and space of the state $u^*$. This is an important factor, since the observation $z=u^*$ will later enter the model reference adaptive system for an update law
\[(q(t),u(t))=MRAS(q_0,z(t)), \qquad q_0: \text{initial guess of } q^*\]
to continuously improve the parameter estimate $q$, which can be viewed as the time continuous version of an approximating sequence to the static parameter $q^*$
and the state estimate $u$ in the sense
\[\|q(t)-q^*\|\to0 \quad\text{and}\quad \|u(t)-z(t)\|\to0  \qquad \text{as } t\to\infty.\]
Constructing an update law undertaking this task, especially when $f$ depends nonlinearly on $q$, is 
in fact the central point in this paper.
\subsubsection*{Notation}
\begin{itemize}[label=]
\item $Q$ is the parameter space, $\cQ$ is the 
time extension of the parameter space. 
\item $U$ is the state space at a fixed time instance,
$\cU$ is the overall state space including time dependence.
\item $Z$ is the observation space at a fixed time instance.
\item $\langle\cdot,\cdot\rangle_{X,X^*}$ denotes the dual paring between $X$ and its dual space. 
\item $\mathcal{L}(X,Y)$ denotes the space of bounded linear operators with operator norm $\|\cdot\|_{X\to Y.}$. For $A\in\mathcal{L}(X,Y)$, we denote by $A^*\in\mathcal{L}(Y^*,X^*)$ its adjoint.
\item $X\hookrightarrow Y$ refers to the continuous embedding with norm $C_{X\to Y}.$
\end{itemize}
}

\subsection{The adaptive system}\label{sec::adptive_sys}
To derive an update law, we first extend the originally stationary
parameter to be a function of time, however, with zero time derivative. This yields the \emph{equivalent model system} 
\begin{alignat}{2}
& D_tq &&= 0  \label{modelq-1}\\
& D_tu^* + f(q,u^*) &&= g  \label{modelq-2}\\
& (q,u^*)(0)= (q^*,u_0). \label{modelq-3}
\end{alignat}
corresponding to the \emph{model system} \eqref{model-1}-\eqref{model-2}, where notationally
we skip the time variable.
Hence, the parameter $q$ behaves as a time constant function that is identical to the exact parameter $q^*$ at each time instance.

Online identification means that the parameter identification, the data collection process and the system operation are taking place simultaneously. During these processes, the data $z$ is fed into a model reference system to adapt its solution $u$.
To this end, we propose the \emph{reference model adaptive} system
\begin{alignat}{2}
& D_tq + \sigma(D_t z+f(q,z)-g)-f'_q(\tilde{q},z)^*(u-z) &&= 0  \label{modelref-1}\\
& D_tu + f(q,z) + \cC(\|q\|_H)(u-z) &&= g  \label{modelref-2}\\
& (q,u)(0)= (q_0,u_0), \label{modelref-3}
\end{alignat}
to be solved for $(q,u).$
\blue{At this point, the parameter estimate $q$ clearly depends on time. The reference adaptive system  is compatible with the equivalent model system in the sense that  compared to \eqref{modelq-1}-\eqref{modelq-2}, the emerging terms in \eqref{modelref-1}-\eqref{modelref-2} continuously decrease the discrepancy $(u(t)-z(t))$. This leads to a progressive improvement of $q(t)$ from $q_0$ towards the exact parameter $q^*$ as times evolves. To fulfill this, $\sigma$ and $\cC$ cannot be arbitrary, but must be chosen appropriately, as the convergence analysis will show.}\\
In the adaptive system, $z$ is the solution to the model system \eqref{model-1}-\eqref{model-2} measured by the observation process \eqref{model-3}.  
$q_0$ is an initial guess for $q^*$, the exact parameter in the model system; \blue{both $q_0$ and $q^*$ are stationary parameters in $Q$. On the other hand, the parameter estimate $q$ lies in some time-dependent parameter space $\cQ$, and $\tilde{q}\in\cQ$ is some point in which $f$ is differentiable. The state estimate $u$ belongs to some state space $\cU$. $\cQ$ and $\cU$ will be detailed later.}\\
\blue{Following usual partial derivative notation, but without strictly demanding Fr\'{e}chet differentiablility,}
we assume 
\begin{align*}
f'_q(\tilde{q},z)\in L^\infty(\I;\cL(H,U^*))
\end{align*}
to be a linearization of $f$ with respect to its first variable; typically just the G\^ateaux derivative. \blue{Here, the inclusion $f'_q(\tilde{q},z)(t)\in\cL(H,U^*)$ corresponds to the bound of $\|q-q^*\|_H$ in $H$ as assumed in \ref{A-coe}.}
Moreover, $\sigma\in\{0,1\}$ is a switching parameter, which is set to zero if $f(q,z)-f(q^*,z)- f'_q(\tilde{q},z)(q-q^*)=0$ (which is, e.g., the case if $f$ is linear with respect to $q$) and to one \blue{otherwise}. \blue{Equally important, $q\mapsto\cC(\|q\|_H)$ is possibly
nonlinear and for fixed $q(t)$,} 
$\cC(\|q(t)\|_H)$ is a linear operator in $\mathcal{L}(U, U^*)$, whose structure will later be specified in Assumption \ref{ass-1}. \blue{The adjoint of the linearization of $f$, the nonlinear term $\cC(\cdot)$ and the design of the function $\sigma(\cdot)$ in our arguments tackle the nonlinearity of the model, especially, in the unknown parameter $q$.}
\\
The model reference adaptive system \eqref{modelref-1}-\eqref{modelref-3} is well-defined due to 
\begin{align*}
&g-D_tz=f(q^*,z)\in L^2(\I;U^*)\subseteq L^2(\I;Q^*)\\
& f'_q(\tilde{q},z):L^2(\I;H)\to L^2(\I;U^*)\\
& f'_q(\tilde{q},z)^*:L^2(\I;U)\to L^2(\I;H)\subset L^2(\I;Q^*).
\end{align*}
This results in $D_tq$ living in the space $L^2(\I;Q^*)$, which is a conventional choice for image spaces in the context of time-dependent PDEs. 

We search for the solution of the adaptive system in the spaces
\begin{alignat}{2}
&u\in \cU:=L^2(\I;U)\cap H^1(\I;U^*)\cap L^\infty(\I;H) \\
&q \in \cQ:=L^2(\I;H)\cap H^1(\I;Q^*)\cap L^\infty(\I;H). 
\end{alignat}
The reference system \eqref{modelref-1}-\eqref{modelref-3} mimics the model system \eqref{modelq-1}-\eqref{modelq-3} in the sense that as time progresses, \eqref{modelref-2} adapts to \eqref{modelq-2} and \eqref{modelref-1} evolves to $D_tq+f(q,z)-f(q^*,z)=0$, driving $q(t)$ towards $q^*$ in \eqref{modelq-1}. Hence, we expect to obtain convergence of the state and of the parameter as $t$ tends to infinity, that is,
\begin{align}\label{convergence}
\errst(t):=u(t)-z(t) \to 0 \quad\text{and}\quad \errpar(t):=q(t)-q^*\to 0 \quad\text{as } t\to\infty.
\end{align}

In order to establish the convergence \eqref{convergence}, we make the following main assumptions.
\begin{assumption}\quad\label{ass-1}
\begin{enumerate}[label=(A\arabic*)]
\item For the exact parameter $q^*\in Q$, the exact state $u^*\in \cU$ \blue{uniquely} exists, and  $z=u^*$ is the full measurement data.\label{A-extU}
\item 
\blue{For $\tilde{q}\in \cQ$ the operator $f'_q(\tilde{q},z)\in L^\infty(\I;\cL(H,U^*))$ is a linearization of $f$ at $\tilde{q}$ in the sense that for any $q\in H$,}
\label{A-lip}
\[
\| f(q,z)-f(q^*,z)- f'_q(\tilde{q},z)(q-q^*)\|_{U^*} \leq L(\|q\|_H)\|q-q^*\|_H.
\] 
\blue{for some monotonically increasing function $L:\I\to\I$.} 
\item If in \ref{A-lip} one has $L\not\equiv0$, i.e., $\sigma=1$, then there exists a constant $\Ccoe>0$ such that for all $q\in Q$, \label{A-coe}
\[\lan f(q,z)-f(q^*,z),q-q^*\ran_{Q^*,Q} \geq \Ccoe\|q-q^*\|_H^2.\] 
\item $\cC$ is chosen such that \label{A-C} for given $q\in H$, and all $v, w\in U$
\label{A-cC}
\begin{align}
\lan \cC(\|q\|_H)v,v \ran_{U^*,U}&\geq\left( \frac{L^2(\|q\|_H)}{2\Ccoe}+M\right)\|v\|_U^2 =:\Mtil(\|q\|_H)\|v\|_U^2
\label{cCCcoeM}\\
\lan \cC(\|q\|_H)v,w \ran_{U^*,U}&\leq\Ntil(\|q\|_H)\|v\|_U\|w\|_U
\end{align}
for some constant $M>0$ and some monotonically increasing function $\Ntil:\I\to\I$.
\end{enumerate}
\end{assumption}

Indeed, the model reference adaptive system \eqref{modelref-1}-\eqref{modelref-3} is an extension of the system devised in \cite{Baumeister-etal} for the case of $f$ being linear with respect to $q$, to the case of nonlinear and possibly also time dependent $f$, that is, of a nonautonomous system (cf. Remark \ref{rem:non-autonomous} below). In our notation, this reads as
\begin{alignat}{2}
& D_tq - f(\cdot,z)^*(u-z) &&= 0  \label{modelref_lin-1}\\
& D_tu + f(q,z) + \cC(u-z) &&= g  \label{modelref_lin-2}\\
& (q,u)(0)= (q_0,u_0). \label{modelref_lin-3}
\end{alignat}
As shown in \cite{Baumeister-etal}, the operator $\cC$ can be chosen as being independent of $q$, which can also be seen as a special case of \ref{A-C}.
In \eqref{modelref-1}-\eqref{modelref-3} as well as in \eqref{modelref_lin-1}-\eqref{modelref_lin-3}, the terms containing $(u-z)$ take into account the residual in the observation equation and exploit it for driving the estimated parameter-state pair $(q,u)$ towards the exact one.
The additional term $D_tz+f(q,z)-g$ in \eqref{modelref-1} is just the residual in the state equation; in case $\sigma=1$ it is used, together with the assumed coercivity \ref{A-coe}, to control the nonlinearity of $f$.

\subsection*{\blue{Discussion of assumptions}}
\begin{remark}
\blue{While existence and uniqueness of the exact state $u^*$ is assumed in \ref{A-extU}, we will prove well-posedness of the model reference adaptive system \eqref{modelref-1}-\eqref{modelref-3} in Section \ref{sec:uniquesolveMRAS} in the general setting. Moreover, \ref{A-extU} is verified for each of the examples in Section \ref{sec:ex} (as are all the other items of Assumption \ref{ass-1}).}
\end{remark}

\begin{remark}
Assumption \ref{A-lip} is satisfied with $L\equiv0$ if $f$ is linear with respect to $q$. In this case, assumption \ref{A-coe} is not needed.
For simplicity of exposition, and since the linear case has already been discussed in \cite{Baumeister-etal}, in the following we only consider the nonlinear case $L\not\equiv0$, $\sigma=1$.
\end{remark}

\begin{remark}\label{remk-strongnorm}
If \ref{A-coe} holds for a stronger norm of $q-q^*$, we also state \ref{A-lip} with this stronger norm on $q-q^*$. In this case, it is more feasible for \ref{A-lip} to be satisfied, thus enabling higher nonlinearity with respect to $q$ in the model.
Moreover, the fact that $z$ belongs to a smoother space similarly enables higher nonlinearity.
\end{remark}

\begin{remark}
The estimate in assumption \ref{A-lip} is fulfilled if, for example, $f$ is locally Lipschitz continuous with respect to $q$, in the sense that
\[\forall M, \exists L(M)\geq0: \|f(q,z)-f(q^*,z)\|_{U^*}\leq L(M)\|q-q^*\|_H \quad \forall q,q^*\in Q: \|q\|_H,\|q^*\|_H\leq M.\]
In addition, $f(\cdot,z)$ is G\^ateaux differentiable on $L^\infty(\I;H)\cap L^2(\I;H)$.\\
Indeed, let us consider, \blue{for any $q$ and $\xi$ with $\|\xi\|\leq1$ in $L^\infty(\I;H)\cap L^2(\I;H)$,}
\begin{equation*}
\blue{R_\epsilon:=}\frac{1}{\epsilon}\| f(q+\epsilon\xi,z) - f(q,z)
- \epsilon f'_q(q,z)\xi \|_{L^2(\T,U^*)}
\blue{=:} \left(\int_0^\infty r_\epsilon(t)^2 dt \right)^\frac{1}{2}.
\end{equation*}
\blue{where $f'_q(q,z)\xi=\lim_{\tau\to0}\frac{1}{\tau}(f(q+\tau\xi,z)-f(q,z))$}
Under the assumption of local Lipschitz continuity, by choosing $M=\|q\|_{L^\infty(\T,H)}+1$, we deduce
\begin{align*}
\|f'_q(q(t),z(t))\xi(t)\|_{U^*}
&=\lim_{\epsilon \rightarrow 0}\left \|\frac{f(q(t)+\epsilon\xi(t),z(t))-f(q(t),z(t))}{\epsilon}\right\|_{U^*} \\
&\leq  L(M)\|\xi(t)\|_H. 
\end{align*}
Thus, for all $\epsilon\in [0,1], \lambda\in [0,1]$,
\begin{align*}
&\|f'_q(q(t)+\lambda\epsilon\xi(t),z(t))\|_{H \rightarrow U^*} \leq L(M).
\end{align*}
\blue{Together with an application of the Mean Value Theorem}, this implies
\begin{align}
r_\epsilon(t)=\left\| \int_0^1 
 (f'_q(q(t)+\lambda\epsilon \xi(t),z(t)) - f'_q(q(t),z(t))) \xi \, d\lambda \right\|_{U^*} \leq  2L(M)\|\xi(t)\|_H,
\end{align}
that is, $r_\epsilon$ is uniformly bounded in $\epsilon$ and square integrable in time. Applying Lebesgue's Dominated Convergence Theorem yields \blue{$R_\epsilon\to 0$ when $\epsilon\to 0$}, thus proving G\^ateaux differentiability of $f(\cdot,z)$ on $L^\infty(\I,H)\cap L^2(\I;H)$.
\end{remark}

\begin{remark}
The coercivity assumption \ref{A-coe}, which we verify on multiple occasions in Section \ref{sec:ex}, considerably facilitates the proof of parameter convergence, as compared to the more general proof via persistence of excitation in \cite{Baumeister-etal}.
To compare these two conditions, we recall that for $f$ linear with respect to $q$, the state $u$ is called uniformly persistently excited iff
\begin{align*}
&\exists \ell>0, \mu>0\, \forall h\in Q\setminus \{0\} \, \forall t_0\in\I\, \exists t_1,t_2\in[t_0,t_0+\ell]\, \exists v\in U\setminus \{0\} \, :\\
&\hspace*{4cm}\left|\int_{t_1}^{t_2}\lan f(h,u(t)),v\ran_{U^*,U}\, ds \right|\geq\mu\|h\|_Q\|v\|_U\,.
\end{align*}
On the other hand, it is readily verified, that by choosing $v:=h\in Q\subseteq U$, $t_1=t_0$, $t_2=t_0+\ell$, and setting  $\mu=\ell \Ccoe$, assumption \ref{A-coe} yields
\begin{align*}
&\exists \ell>0, \mu>0\, \forall h\in Q\setminus \{0\} \, \forall t_0\in\I\, \exists t_1,t_2\in[t_0,t_0+\ell]\, \exists v\in U\setminus \{0\} \, :\\
&\hspace*{4cm}\left|\int_{t_1}^{t_2}\lan f(h,u(t)),v\ran_{U^*,U}\, ds \right|\geq\mu\|h\|_H\|v\|_H\,.
\end{align*}
\ref{A-coe} indeed implies a certain persistence of excitation, albeit with respect to weaker norms. 
\end{remark}

\begin{remark} \label{rem:non-autonomous}
The proposed method can be generalized to nonautonomous model systems with $f:\I\times H\times U \to U^*$. In this case, we assume that $f$ satisfies the Carath\'eodory conditions and \ref{A-lip}, \ref{A-coe} hold uniformly for almost all $t\in(0,\infty).$  
\end{remark}

\section{Convergence analysis of the error system}\label{sec:Error}
By subtracting the model system \eqref{modelq-1}-\eqref{modelq-3} from the reference system \eqref{modelref-1}-\eqref{modelref-3}, we see that the error components $(\errst,\errpar)=(u-u^*,q-q^*)$ satisfy the nonlinear \emph{error system}
\begin{alignat}{2}
&D_t\errst + f(q,z)-f(q^*,z)+\cC(\|q\|_H)\errst \quad&&=0 \label{modelerr-1}\\
&D_t\errpar + f(q,z)-f(q^*,z)-f'_q(\tilde{q},z)^*\errst &&=0  \label{modelerr-2}\\
&(\errst,\errpar)(0)= (0,q_0-q^*)  \label{modelerr-3}
\end{alignat}

In the following section, we establish the convergence of $(\errst,\errpar)$ \blue{in the noise-free case.}
\subsection{Convergence of $(\errst,\errpar)$}
We test \eqref{modelerr-1} and \eqref{modelerr-2} respectively by $\errst(t)\in U$ and $\errpar(t)\in Q$, then sum up the outcome
\begin{align}
0&=\frac{1}{2}\frac{d}{dt}\left[ \|\errst\|_H^2+\|\errpar\|_H^2 \right](t) + \lan f(q(t),z(t))-f(q^*(t),z(t)),\errpar(t) \ran_{Q^*,Q} \nonumber\\
&\hspace{1cm}+ \lan f(q(t),z(t))-f(q^*(t),z(t))-f'_q(\tilde{q},z(t))\errpar(t),\errst(t) \ran_{U^*,U} + \lan \cC(\|q(t)\|_H)\errst(t),\errst(t) \ran_{U^*,U} \nonumber\\
&\geq \frac{1}{2}\frac{d}{dt}\left[ \|\errst(t)\|_H^2+\|\errpar(t)\|_H^2 \right] + \Ccoe\|\errpar(t)\|_H^2 - L(\|q(t)\|_H)\|\errpar(t)\|_H\|\errst(t)\|_U \nonumber\\
&\hspace*{5cm}+ \Mtil(\|q(t)\|_H)\|\errst(t)\|_U^2 \nonumber\\
&\geq \frac{1}{2}\frac{d}{dt}\left[ \|\errst(t)\|_H^2+\|\errpar(t)\|_H^2 \right] + \Ccoe\|\errpar(t)\|_H^2 - \left(\frac{\Ccoe}{2}\|\errpar(t)\|_H^2 +\frac{L^2(\|q(t)\|_H)}{2\Ccoe}\|\errst(t)\|_U^2 \right)  \nonumber\\
&\hspace*{5cm}+ \Mtil(\|q(t)\|_H)\|\errst(t)\|_U^2 \nonumber\\
&\geq \frac{1}{2}\frac{d}{dt}\left[ \|\errst\|_H^2+\|\errpar\|_H^2 \right](t) + \frac{\Ccoe}{2}\|\errpar(t)\|^2_H  + M\|\errst(t)\|_U^2. \label{bound}
\end{align}
Above, we make use of assumptions \ref{A-coe}-\ref{A-C}. This estimate reveals the first observation 
\begin{align*}
t\mapsto  E(t):=\|\errst(t)\|_H^2+\|\errpar(t)\|_H^2 \quad\text{is decreasing,}
\end{align*}
on the error dynamics, which implies the state- and the parameter-error stay bounded.\\
Furthermore, we have the inequality
\[E'(t) + \min\left\{\Ccoe;2MC_{U\to H}\right\}E(t)\leq 0,\]
where $C_{U\to H}$ is the norm of the continuous embedding $U\hookrightarrow H$. This means that the function defined by $\tilde{E}(t):=\exp\big(\min\left\{\Ccoe;2MC_{U\to H}\right\}\,t\big)E(t)$ satisfies $\tilde{E}'(t)\leq0$ and is therefore monotonically decreasing. In particular, $\tilde{E}(t)\leq \tilde{E}(0)=\left[ \|\errst(0)\|_H^2+\|\errpar(0)\|_H^2 \right]$ implies
\begin{align}\label{rate_exp}
E(t)\leq \exp\big(-\min\left\{\Ccoe;2MC_{U\to H}\right\}\,t\big)\left[ \|\errst(0)\|_H^2+\|\errpar(0)\|_H^2 \right],
\end{align}
that is, an exponential convergence rate of the total error.

Next, by integrating \eqref{bound} over $[0,T]$ we get
\begin{align}\label{bound-all}
&\frac{1}{2}\left[ \|\errst(T)\|_H^2+\|\errpar(T)\|_H^2 \right]+\frac{\Ccoe}{2}\int_0^T\|\errpar(t)\|^2_H\,dt + M \int_0^T\|\errst(t)\|_U^2\,dt \nonumber\\
&\hspace{1cm} \leq \frac{1}{2}\left[ \|\errst(0)\|_H^2+\|\errpar(0)\|_H^2 \right] \qquad \text{for any }T>0
\end{align}
thus
\begin{align}\label{bound-et}
&\|D_t\errst\|_{L^2(\I,U^*)}^2=\int_0^\infty\left(\sup_{\|v\|_U\leq1}\lan D_t\errst , v \ran_{U^*,U}\right)^2\, dt \nonumber\\
&=\int_0^\infty\left(\sup_{\|v\|_U\leq1}\lan -f(q,z)+f(q^*,z)-\cC(\|q\|_H)\errst , v \ran_{U^*,U}\right)^2\, dt \nonumber\\
&\leq \int_0^\infty \Bigl(\sup_{\|v\|_U\leq1}\Bigl[ \Bigl(L(\|q(t)\|_H) + \|f'_q(\tilde{q},z(t))\|_{H\to U^*}\Bigr)\|\errpar(t)\|_H\|v\|_U \nonumber\\
&\hspace*{6cm}+ \Ntil(\|q(t)\|_H)\|\errst(t)\|_U\|v\|_U \Bigr]\Bigr)^2\,dt \nonumber\\
&\leq \frac{(L(\|q\|_{L^\infty(\I;H)})+\|f'_q(\tilde{q},z)\|_{L^\infty(\I;\cL(H,U^*))} + \Ntil(\|q\|_{L^\infty(\I;H)}))^2}{\min\left\{\Ccoe;2M\right\}}
 \nonumber\\
&\hspace*{6cm}
\times \left[ \|\errst(0)\|_H^2+\|\errpar(0)\|_H^2 \right.]
\end{align}
Additionally, due to $Q\subseteq U$, it follows that
\begin{align}\label{bound-rt}
&\|D_t\errpar\|_{L^2(\I;Q^*)}\leq \|D_t\errpar\|_{L^2(\T,U^*)}=\int_0^\infty\sup_{\|p\|_U\leq1}\lan D_t\errpar , p \ran_{U^*,U}\, dt \nonumber\\
&=\int_0^\infty\left(\sup_{\|p\|_U\leq1}\lan -f(q,z)+f(q^*,z)+f'_q(\tilde{q},z)^*\errst , p \ran_{U^*,U}\right)^2\, dt \nonumber\\
&\leq \int_0^\infty  \Bigl(\sup_{\|p\|_U\leq1}\Bigl[ (L(\|q(t)\|_H) + \|f'_q(\tilde{q},z)(t)\|_{H\to U^*})\|\errpar(t)\|_H\|p\|_U \nonumber\\
&\hspace*{6cm}+ \|f'_q(\tilde{q},z)(t)\|_{H\to U^*}\|p\|_H\|\errst(t)\|_U \Bigr]\Bigr)^2\,dt \nonumber\\
&\leq \frac{(L(\|q\|_{L^\infty(\I;H)})+2\|f'_q(\tilde{q},z)\|_{L^\infty(\I;\cL(H,U^*))} )^2}{\min\left\{\Ccoe;2M\right\}}\,\left[ \|\errst(0)\|_H^2+\|\errpar(0)\|_H^2 \right]\,.
\end{align}

Summarizing  \eqref{rate_exp}-\eqref{bound-rt}, we state the result:
\begin{proposition}\label{prop-noisefree}
Let Assumption \ref{ass-1} be fulfilled. Then the following statements on the parameter $q$ and the state $u$ as well as the corresponding errors $\errpar=q-q^*$, $\errst=u-u^*$ hold true:,
\begin{enumerate}[label=(\roman*)]
\item
$u\in \cU=L^2(\I;U)\cap H^1(\I;U^*)\cap L^\infty(\I;H),$ \\
$q \in \cQ=L^2(\I;H)\cap H^1(\I;Q^*)\cap L^\infty(\I;H).$
\item
\begin{align*}
&\sup_{t\geq 0}\left[ \|\errst(t)\|_H^2+\|\errpar(t)\|_H^2 \right]+\Ccoe\int_0^\infty\|\errpar(s)\|^2_H\,ds + 2M\int_0^\infty\|\errst(s)\|_U^2\,ds \nonumber\\
&\hspace{1cm} \leq \left[ \|\errst(0)\|_H^2+\|\errpar(0)\|_H^2 \right]. 
\end{align*}
\item For all $t\geq0$,
\begin{align*}
\left[ \|\errst(t)\|_H^2+\|\errpar(t)\|_H^2 \right]\leq \exp\big(-\min\left\{\Ccoe;2MC_{U\to H}\right\}t\big)\left[ \|\errst(0)\|_H^2+\|\errpar(0)\|_H^2 \right].
\end{align*}
\end{enumerate}
\end{proposition}

\def\deltil{\tilde{\delta}}
\def\ti{{\rm ti}}
\def\sp{{\rm sp}}
\def\cR{\mathcal{R}}

\subsection{Case of noisy data}
We now turn to the practically relevant setting of noisy observations. In other words, we have $z^\delta$ in place of $z=u^*$, with a certain noise level $\delta>0$ that we assume to be given in the data space norm, that is, 
\begin{equation}\label{delta}
\|z^\delta-z\|_{L^p(\I;Z)}\leq \delta\,.
\end{equation}
Since $z^\delta$ does not satisfy the regularity requirements needed to be inserted into the model reference adaptive system, we smooth it by, e.g. filtering or local averaging, or, more abstractly, by applying regularizing operators $\cR^\sp:Z\to U$ (pointwise in time), $\cR^\ti:L^p(\I;Z)\to W^{1,p}(\I;H)$ (typically nonlocal in time) such that -- by an appropriate choice of the regularization parameters contained in the definition of $\cR^\sp$, $\cR^\ti$ -- the estimates 
\begin{equation}\label{deltatil}
\|\cR^\sp(z^\delta(t))-z(t)\|_{U}\leq \deltil^\sp(t) \,, \quad 
\|D_t(\cR^\ti(z^\delta)-z)\|_{L^p(\I;H)}\leq \deltil^\ti
\end{equation}
hold for some $p\geq2$. Inserting the smoothed data in place of $z$, we get, instead of \eqref{modelref-1}-\eqref{modelref-3},
\begin{alignat}{2}
& D_tq + \sigma(D_t \cR^\ti(z^\delta)+f(q,\cR^\sp(z^\delta))-g)-f'_q(\tilde{q},\cR^\sp(z^\delta))^*(u-\cR^\sp(z^\delta)) &&= 0  
\nonumber\\
& D_tu + f(q,\cR^\sp(z^\delta)) + \cC(\|q\|_H)(u-\cR^\sp(z^\delta)) &&= g  
\nonumber\\
& (q,u)(0)= (q_0,u_0)
\nonumber\,.
\end{alignat}
Again, we focus on the case $\sigma=1$. In addition to \eqref{deltatil}, we make the following assumptions.
\begin{assumption}\quad\label{ass-2}
\begin{enumerate}[label=\blue{(B\arabic*)}]
\item There exists a constant $\tilde{L}_0$ such that for all $v,w\in U,$
\[
\|f'_q(\tilde{q},v)-f'_q(\tilde{q},w)\|_{\cL(H,U^*)}\leq \tilde{L}_0 \|v-w\|_U\,.
\]
\item There exist constants $\tilde{L}_1$, $\tilde{L}_2$ such that for all $q\in H$,  $v,w\in U,$
\[
\|f(q,v)-f(q,w)\|_{U^*}\leq \tilde{L}_1 \|v-w\|_U\,, \quad 
\|f(q,v)-f(q,w)\|_{H}\leq \tilde{L}_2 \|v-w\|_U\,.
\]
\end{enumerate}
\end{assumption}
Moreover, we replace the condition \eqref{cCCcoeM} on the choice of $\cC$ in \blue{Assumption \ref{ass-1}} \ref{A-cC} by 
\begin{align}
\lan \cC(\|q\|_H,t)v,v \ran_{U^*,U}&\geq\left( \frac{(L(\|q\|_H)+\tilde{L}_0\|\deltil^\sp\|_{L^\infty(0,t)})^2}{2\Ccoe}+M\right)\|v\|_U^2 =:\Mtil(\|q\|_H,t)\|v\|_U^2\,.
\label{cCCcoeM_noi}
\end{align}
\begin{proposition}\label{prop-noisy}
Let Assumption \ref{ass-1} with \eqref{cCCcoeM_noi} in place of \eqref{cCCcoeM}, Assumption \ref{ass-2}, and the noise bound \eqref{deltatil} be satisfied. Then the following statements on the parameter $q$ and the state $u$ as well as the corresponding errors $\errpar=q-q^*$, $\errst=u-u^*$ hold true:
\begin{enumerate}[label=(\roman*)]
\item
$u\in \cU=L^2(\I;U)\cap H^1(\I;U^*)\cap L^\infty(\I;H),$ \\
$q \in \cQ=L^2(\I;H)\cap H^1(\I;Q^*)\cap L^\infty(\I;H).$
\item For any $\omega < \min\{\Ccoe,2MC_{U\to H}\},$ there exists $C>0$ such that for all $t\geq0$,
\begin{align*}
\left[ \|\errst(t)\|_H^2+\|\errpar(t)\|_H^2 \right]\leq \exp\big(-\omega t\big)\left[ \|\errst(0)\|_H^2+\|\errpar(0)\|_H^2 \right]
+ C \Bigl(\|\deltil^\sp\|_{L^p(0,t)}^2 + (\deltil^\ti)^2\Bigr).
\end{align*}
\end{enumerate}
\end{proposition}
\proof 
The crucial estimate follows analogously to the proof of Proposition \ref{prop-noisefree}, using the fact that with noisy data, the error system becomes 
\begin{alignat}{2}
&D_t\errst + f(q,z)-f(q^*,z)+\cC(\|q\|_H)\errst \quad&&=d^u \label{modelerr-1_noi}\\
&D_t\errpar + f(q,z)-f(q^*,z)-(f'_q(\tilde{q},z) - d^0)^*\errst 
&&=d^q  \label{modelerr-2_noi}\\
&(\errst,\errpar)(0)= (0,q_0-q^*)  \label{modelerr-3_noi}
\end{alignat}
in place of \eqref{modelerr-1}-\eqref{modelerr-3}, where 
\[
\begin{aligned}
d^0(t)=& f'_q(\tilde{q},z(t))-f'_q(\tilde{q},\cR^\sp(z^\delta(t)))\, \in \cL(H,U^*)\\
d^u(t)=& f(q(t),z(t))-f(q(t),\cR^\sp(z^\delta(t))) - \cC(\|q(t)\|_H)(z(t)-\cR^\sp(z^\delta(t)))\ \in U^*\\
d^q(t)=& f(q(t),z(t))-f(q(t),\cR^\sp(z^\delta(t))) + D_t(z(t)-\cR^\ti(z^\delta(t)))\\
&\hspace*{4cm}+ f'_q(\tilde{q},\cR^\sp(z^\delta(t)))^*(z(t)-\cR^\sp(z^\delta(t)))\ \in H.
\end{aligned}
\] 
By testing \blue{\eqref{modelerr-1_noi} and \eqref{modelerr-2_noi}} with $\errst(t)$ and $\errpar(t)$, respectively, we get, in place of \eqref{bound}, 
\begin{align}
&\lan d^u(t),\errst(t)\ran_{U^*,U}+\lan d^q(t),\errpar(t)\ran_{Q^*,Q}\nonumber\\
&\geq \frac{1}{2}\frac{d}{dt}\left[ \|e\|_H^2+\|r\|_H^2 \right](t) + \Ccoe\|\errpar(t)\|_H^2\nonumber\\ 
&\qquad- \Bigl(L(\|q(t)\|_H+ \tilde{L}_0\deltil^\sp(t)\Bigr)\|\errpar(t)\|_H \|\errst(t)\|_U 
+ \Mtil(\|q(t)\|_H)\|\errst(t)\|_U^2 \nonumber\\
&\geq \frac{1}{2}\frac{d}{dt}\left[ \|\errst(t)\|_H^2+\|\errpar(t)\|_H^2 \right] + \Ccoe\|\errpar(t)\|_H^2 
\nonumber\\
&\qquad- \left(\frac{\Ccoe}{2}\|\errpar(t)\|_H^2 +\frac{(L(\|q(t)\|_H)+ \tilde{L}_0\deltil^\sp(t))^2}{2\Ccoe}\|\errst(t)\|_U^2 \right)  + \Mtil(\|q(t)\|_H) \|\errst(t)\|_U^2 
\nonumber\\
&\geq \frac{1}{2}\frac{d}{dt}\left[ \|\errst(t)\|_H^2+\|\errpar(t)\|_H^2 \right] + \frac{\Ccoe}{2}\|\errpar(t)\|^2_H  + 
M\|\errst(t)\|_U^2. \label{bound_noi}
\end{align}
An application of Young's Inequality and Assumption~\ref{ass-2} as well as multiplication by two, it yields
\[
\begin{aligned}
E'(t)+\omega E(t) 
&\leq 
 \frac{1}{\epsilon_1} \|d^u(t)\|_{U^*}^2 + \frac{1}{\epsilon_2} \|d^q(t)\|_{H}^2\\
&\leq C(\epsilon_0,\epsilon_1,\epsilon_2) \deltil^\sp(t)^2+ \frac{2}{\epsilon_2}\|D_t(z(t)-\cR^\ti(z^\delta)(t))\|_{H}^2 =: D(t)
\,,
\end{aligned}
\]
where $\omega=\min\left\{\Ccoe-\epsilon_2;2MC_{U\to H}-\epsilon_1\right\}$, 
\[
C(\epsilon_0,\epsilon_1,\epsilon_2)=\sup_{t\in\I}\frac{(\tilde{L}_1+\Ntil(\|q(t)\|_H))^2}{\epsilon_1}  + \frac{2}{\epsilon_2} \bigl(\tilde{L}_2 + \|f'_q(\tilde{q},z(t))\|_{\cL(H,U^*)}+\tilde{L}_0\deltil^\sp(t)\bigr)^2\,.
\]
For $\tilde{E}(t):=e^{\omega t} E(t),$ this means $\tilde{E}'(t)\leq e^{\omega t} D(t)$ and thus 
$\tilde{E}(t)\leq \tilde{E}(0)+\int_0^t e^{\omega s} D(s)\, ds=E(0)+\int_0^t e^{\omega s} D(s)\, ds$, that is, 
\[
\begin{aligned}
E(t)&\leq e^{-\omega t} E(0)+\int_0^t e^{-\omega (t-s)} D(s)\, ds
\leq e^{-\omega t} E(0)+ C_\omega \|D\|_{L^{p/2}(0,t)} \\
&\leq e^{-\omega t} E(0)+ C_\omega \Bigl(C(\epsilon_0,\epsilon_1,\epsilon_2) \|\deltil^\sp\|_{L^p(0,t)}^2
+ \frac{2}{\epsilon_2} (\deltil^\ti)^2\Bigr),
\end{aligned}
\]
where $C_\omega = 1$ if $p=2$, and $C_\omega = (\frac{p}{p-2}\omega)^{-\frac{p-2}{p}}$ if $p>2$.
\qed

\subsection{Unique solvability of the adaptive system}\label{sec:uniquesolveMRAS}
To examine unique solvability of the proposed reference adaptive system, we first reformulate it as the initial value problem
\begin{align}
&D_t(u,q)(t)+\F(u,q)(t)=\G(t) \qquad t>0\\
& (u,q)(0)= (u_0,q_0),
\end{align}
where
\begin{align}
&\F:U\times H\to U^*\times Q^* \nonumber\\
&\F:
\begin{pmatrix}
u(t) \\ q(t)
\end{pmatrix}
\mapsto
\begin{pmatrix}
 f(q,z)(t) + \cC(\|q(t)\|_H)(u-z)(t)\\
 f(q,z)(t)-f(q^*,z)(t)-f'_q(\tilde{q},z)^*(u-z)(t)
\end{pmatrix}
=:
\begin{pmatrix}
 \F_1(u(t),q(t))\\
 \F_2(u(t),q(t))
\end{pmatrix},\\[1ex]
&\G(t):=
\begin{pmatrix}
g(t)\\
0
\end{pmatrix}\in U^*\times Q^*.
\end{align}
We then define the pairing
\begin{align*}
\lan \F(u,q),(\util,q^*)\ran:= \lan \F_1(u,q),\util\ran_{U^*,U} + \lan \F_2(u,q),q^*\ran_{Q^*,Q}
\end{align*}
between $U^*\times Q^*$ and $U\times Q.$ \\In the following, we prepare some evaluations, which will be used to prove pseudomonotonicity of $\F$ \cite[Definition 2.1]{Roubicek}.

\subsubsection*{Pseudomonotonicity with respect to $u$}
Consider the function $A:=\cC(\|\cdot\|_H)(\cdot-z)$. Since $A$ is bounded and demicontinuous on  $U\times H$ , i.e. $(u_n,q_n)\to(u,q)$, one has $A(u_n,q_n) \rightharpoonup A(u,q), n\to\infty$. Referring to \cite[-Lemma 2.10]{Roubicek}, $A$ is pseudomonotone if the following statement holds:
\begin{align}
&\text{If}\quad (u_n,q_n) \xrightharpoonup{n\to\infty} (u,q) \quad\text{and}\quad \limsup_{n\to\infty}\,\lan A(u_n,q_n)-A(u,q),(u_n,q_n)-(u,q) \ran\leq 0, \label{pasudo-if}\\
&\text{then}\quad (u_n,q_n) \xrightarrow{n\to\infty}(u,q). \label{pasudo-then}
\end{align}

It suffices to consider the simple form $\cC(\|\cdot\|_H)(\cdot-z)=\|\cdot\|_H\cC(\cdot-z)=:A(e,q)$, the argument for the higher order of $\cC(\|\cdot\|_H)$ being similar. \blue{With \ref{A-C}, we estimate}
\begin{align*}
&\lan A(u_n,q_n)-A(u,q),(u_n,q_n)-(u,q) \ran_{U^*,U}=\lan A(e_n,q_n)-A(e,q),e_n-e \ran_{U^*,U}\\
&=\lan \|q_n\|_H\cC(e_n)-\|q\|_H\cC(e),e_n-e \ran_{U^*,U}\\
&=\|q_n\|_H\|e_n\|^2_U + \|q\|_H\|e\|^2_U 
- \|q_n\|_H\lan \cC(e_n),e\ran_{U^*,U} - \|q\|_H\lan \cC(e),e_n\ran_{U^*,U}\\
&\geq \|q_n\|_H\|e_n\|^2_U + \|q\|_H\|e\|^2_U 
-\frac{\|q_n\|_H}{2}\left( \|e_n\|_U^2+\|e\|_U^2 \right)
-\frac{\|q\|_H}{2}\left( \|e_n\|_U^2+\|e\|_U^2 \right)\\
&=\frac{1}{2} \left( \|q_n\|_H-\|q\|_H \right)\left( \|e_n\|_U^2-\|e\|_U^2 \right).
\end{align*}
The assumption \eqref{pasudo-if} in combination with weak lower semicontinuity of $\|\cdot\|_H$ and $\|\cdot\|_U$ leads to
\begin{align}
0&\geq \limsup_{n\to\infty}\,\lan A(u_n,q_n)-A(u,q),(u_n,q_n)-(u,q) \ran_{U^*,U}\nonumber\\
&\geq \limsup_{n\to\infty}\,\frac{1}{2} \left( \|q_n\|_H-\|q\|_H \right)\left( \|e_n\|_U^2-\|e\|_U^2 \right) \nonumber\\
&\geq \liminf_{n\to\infty}\,\frac{1}{2} \left( \|q_n\|_H-\|q\|_H \right)\left( \|e_n\|_U^2-\|e\|_U^2 \right) \nonumber\\
&=0, \label{pseu-wsc}
\end{align}
which implies the limit $\left( \|q_n\|_H-\|q\|_H \right)\left( \|e_n\|_U^2-\|e\|_U^2 \right)\to0, n\to\infty.$ This shows the fact that as $n\to\infty$, either $e_n\xrightarrow{U} e$, 
or
$q_n\xrightarrow{H} q$, 
\blue{or both}.

In addition, we observe
\begin{align*}
0\xleftarrow{n\to\infty}&= \lan \|q_n\|_H\cC(e_n)-\|q\|_H\cC(e),e_n-e \ran_{U^*,U}\\
&=\lan \|q_n\|_H\cC(e_n)-\|q_n\|_H\cC(e)+\|q_n\|_H\cC(e)-\|q\|_H\cC(e),e_n-e \ran_{U^*,U}\\
&=\|q_n\|_H\|e_n-e\|_U^2 + (\|q_n\|_H-\|q\|_H)\lan e,e_n-e\ran_{U^*,U}.
\end{align*}
Then, if $q_n\xrightarrow{H} q$, the second term in this sum converges to $0$. This constrains $\|e_n-e\|_U$ to converge to $0$ as well. As a consequence,
\begin{align}\label{pseu-u}
\text{As} n\to\infty:\quad e_n \xrightarrow{U} e \quad\text{always.}
\end{align}
Comparing to \eqref{pasudo-then}, $A$ is pseudomonotone with respect to the variable $u.$
\subsubsection*{Pseudomonotonicity with respect to $q$}
Strong convergence of $q_n$ in $H$ is straightforwardly attainable by using \ref{A-coe}:
\begin{align*}
0\geq\limsup_{n\to\infty}\lan f(q_n,z)-f(q,z),q_n-q \ran_{U^*,U} 
\geq \Ccoe\lim_{n\to\infty}\|q_n-q\|^2_H,
\end{align*}
which shows
\begin{align}\label{pseu-q}
\text{As} n\to\infty:\quad q_n \xrightarrow{H} q.
\end{align}
Together with the fact that $f(\cdot,z)$ is bounded and continuous by \ref{A-lip}, pseudomonotonicity of $f$ with respect to the variable $q$ is confirmed.\\

Combining the obtained results, we are now able to evaluate pseudomonotonicity of $\F.$
\subsubsection*{Pseudomonotonicity with respect to $(u,q)$}
First of all, boundedness and continuity of $\F$ is deduced from these properties of $f(\cdot,z)$ and $A$. Then, under the premise \eqref{pasudo-if},
\begin{align*}
&(u_n,q_n) \xrightharpoonup{n\to\infty} (u,q)\\
&\limsup_{n\to\infty}\,\lan \F(u_n,q_n)-\F(u,q),(u_n,q_n)-(u,q) \ran\leq 0,
\end{align*}
and we obtain
\begin{align*}
&K:= \lan \F_1(u_n,q_n)-\F_1(u,q),e_n-e\ran_{U^*,U} + \lan \F_2(u_n,q_n)-\F_2(u,q),q_n-q\ran_{Q^*,Q}\\
&= \lan f(q_n,z)-f(q,z),q_n-q\ran_{Q^*,Q}
+ \lan f(q_n,z)-f(q,z)- f'_q(\tilde{q},z)(q_n-q), e_n-e \ran_{U^*,U}\\
&\hspace{1cm} + \lan \cC(\|q_n\|)e_n-\cC(\|q\|)e, e_n-e\ran_{U^*,U}\\
&=: K_1+K_2+K_3 \\
&\geq \frac{1}{2} K_1+\frac{1}{2} K_3
\end{align*}
by estimating  in the same manner as in \eqref{bound}. From the facts $K_1\geq0$ (see \ref{A-coe}) and $K_3\geq0$ for sufficiently large $n$ (see \eqref{pseu-wsc}), we get the implication
\begin{align*}
\limsup_{n\to\infty} K\leq 0 \quad\text{implies}\quad 
\begin{cases}
\limsup_{n\to\infty} K_1=0 \quad(\text{or}\leq0)\\
\limsup_{n\to\infty} K_3=0 \quad(\text{or}\leq0)
\end{cases}.
\end{align*}
This enables us to apply the achieved result \eqref{pseu-q} for $K_1$ and \eqref{pseu-u} for $K_3$, to thus conclude
\begin{align}\label{pseu-uq}
\text{As} n\to\infty:\quad (u_n,q_n) \xrightarrow{U\times H} (u,q).
\end{align}
With that, pseudomonotonicity of $\F$ with respect to $(u,q)$ has been proven.\\

We are now ready to show unique existence of a solution to the adaptive system \eqref{modelref-1}-\eqref{modelref-2} through the three following steps:
\begin{enumerate}
\item \emph{Approximate solutions}\\
Using the Galerkin method, we can construct approximate solutions to \eqref{modelref-1}-\eqref{modelref-2} on finite-dimensional subspaces.\\
Referring to Proposition \ref{prop-noisefree}, these approximate solutions are uniformly bounded. Thus, there exists a subsequence that weakly converges to some $(u,q)\in\cU\times\cQ.$

\item \emph{Weak limit of approximate solutions is a solution}\\
This weak limit $(u,q)$ is indeed a weak solution to \eqref{modelref-1}-\eqref{modelref-2}, since $\F$ is pseudomonotone by \eqref{pseu-uq}.\\[2ex]
\emph{Proof.} \cite[Theorem 8.27]{Roubicek}

\item \emph{Uniqueness of the solution}\\
\eqref{bound} shows that $F$ satisfies the condition 
\[
\begin{aligned}
&\exists\rho_1, \rho_2\in\R, \forall(u,q), (\util,q^*)\in U\times H:\,\\
&\hspace*{2cm} \lan \F(u,q)-\F(\util,q^*),(u,q)-(\util,q^*)\ran \geq \rho_1\|u-\util\|_H^2+\rho_2\|q-q^*\|_H^2, 
\end{aligned}
\]
which ensures uniqueness of the solution to the adaptive system in $\cU\times \cQ$. \\[2ex]
\emph{Proof.} \cite[Theorem 8.31]{Roubicek}
\end{enumerate}

\section{Examples}  \label{sec:ex}
\blue{Motivated by benchmark PDEs examples, we investigate our proposed online method by means of
the task of determining unknown coefficients in nonlinear parabolic equations. We stay within the framework of the unknown coefficients being spatially (but not temporally)
variant, and full measurements of the exact states associated with exact coefficients are available.} \\
To begin, we analyze unique existence for the perturbed linear parabolic problem  
\begin{alignat}{3}
& D_tu - \nabla\cdot(a\nabla u) + cu + \phi(u)\psi(a,c) = g  \quad&&\mbox{ in }\Omega\times\T\label{ex0-1}\\
& u(0) = u_0 &&\mbox{ in }\Omega\times\{0\},\label{ex0-2}
\end{alignat}
with nonlinear terms $\phi, \psi$. Here $\Omega\subset\R^3$ is a bounded smooth domain, and
\begin{equation}\label{ex0-general}
\begin{aligned}
& c\in L^2(\Omega),\,\,\, \qquad 0<\underline{a}\leq a(x)\leq \overline{a}\quad \forall^\text{a.e.} x\in \Omega\,\
&U=H^1(\Omega),\qquad H=L^2(\Omega,)\\
& g\in W^{1,\infty,2}(\T;U^*,U^*),\\
& u_0\in U.
\end{aligned}
\end{equation}

In the following, we present the unique existence result.
\begin{proposition}\label{prop-assum}
In Problem \eqref{ex0-1}-\eqref{ex0-2}, we assume
\begin{align*}
&a\in L^\infty(\Omega)\cap W^{1,3}(\Omega),\quad a>0 \text{ a.e. on }\Omega, \\
& c\in L^2(\Omega), \quad\text{with}\quad C_{H^1\to L^6}\|\nabla a\|_{L^3(\Omega)} + C_{H^2\to L^\infty}\|\tilde{c}-c\|_{L^2(\Omega)} < \underline{a}/C^{\tilde{c},a} \quad\text{for some }\tilde{c}\in L^\infty(\Omega),\\
& g\in W^{1,\infty,2}(\T;H,U^*),\\
& u_0 \in H^2(\Omega),\\
& \psi(c, a)\in L^\infty(\Omega),\\
& \phi(\cdot)\psi(a,c)\, \text{ is monotone,}\qquad\text{and}\qquad |\phi(u)|\leq C_\phi(1+|u|^3).
\end{align*}
Then the initial value problem \eqref{ex0-1}-\eqref{ex0-2} admits a unique solution 
\[u\in \cU_\infty:= W^{1,\infty,\infty}(\T;H,H)\cap W^{1,\infty,2}(\T;U,U)\cap L^\infty(\T\times\Omega)\quad \forall T>0.\]
\end{proposition}
\proof
\blue{\nameref{appendix}}
\qed\\

By default, these conditions will be imposed on the model equations in the upcoming examples. This is to ensure that the exact state, which is as always denoted by $z$, exists in the right space as stated in \ref{A-extU}. Moreover, its uniform boundedness in $(t,x)$ facilitates the validation of assumptions \ref{A-lip}-\ref{A-coe} in the sense that positivity (or negativity) of $z$ features there. With respect to this existence result and regularity of $z$, we assume
\begin{align}\label{ex0-boundz}
\text{for any } T>0, \exists M_z>0, \forall(x,t)\in\Omega\times(0,T): \quad|\phi(z)|\leq C_\phi(1+|z|^3)\leq M_z.
\end{align}
The coefficients $c, a$ here will respectively play the role of the exact space-dependent coefficient $q^*$ in the examined problems. Thus, we can impose
\begin{align}\label{ex0-boundq}
\exists M_{q^*}>0, \forall x\in\Omega: \quad|\psi(q^*)|\leq M_{q^*}.
\end{align}

\subsection{Identification of a potential -- $c$ problem with nonlinear perturbation}
We consider the estimation of the parameter $q$ in the model 
\begin{align}
& D_tu - \Delta u + qu + \phi(u)\psi(q) = g  \quad\mbox{ in }\Omega\times\I\label{ex1-1}\\
& u(0) = u_0 \mbox{ in }\Omega\times\{0\}\label{ex1-2}
\end{align}
on a smooth bounded domain $\Omega\subseteq\mathbb{R}^3$ under the assumptions
\begin{align*}
&\text{the exact state } z(t,x)\geq\underline{c}>0  \quad \forall x\in\Omega,\,t>0,\\
&\textcolor{black}{\phi(z)\psi(\cdot)} \text{ is monotone}, \quad\text{and}\quad
|\psi'(q)|\leq C_\psi(1+|q|^{\beta-1}).
\end{align*}
To achieve boundedness away from zero of $z$ by means of a maximum principle \textcolor{black}{ \cite[Lemma 2.1, Chapter 2]{Pao}}, we impose inhomogeneous Dirichlet boundary conditions
\[
z=h\geq\underline{c} \mbox{ on }\partial\Omega\times\I
\]
and assume that also 
\[
u_0\geq\underline{c}\,, \quad g\geq q^*\underline{c}\textcolor{black}{\, + M_zM_{q^*}}\,.
\]
From this and the fact that $\hat{z}=z-\underline{c}$ solves
\begin{alignat*}{3}
& D_t\hat{z} - \Delta \hat{z} + q^*\hat{z} = g-q^*\underline{c}-\phi(z)\psi(q^*)\geq0  \quad&&\mbox{ in }\Omega\times\T\\
&\hat{z}=h-\underline{c}\geq0 &&\mbox{ on }\partial\Omega\times\I\\
&\hat{z}(0) = u_0-\underline{c}\geq0 &&\mbox{ in }\Omega\times\{0\},
\end{alignat*}
we conclude $\hat{z}\geq0$, i.e., $z\geq\underline{c}$. \textcolor{black}{Since $z\in L^\infty(\T\times\Omega), \forall T>0$ as claimed in Proposition \ref{prop-assum}, it makes sense to have $z\geq\underline{c}$ on $\T, \forall T$; hence, $z\geq\underline{c}$ on $\I$.}
\\
In order to work with homogeneous boundary conditions, we now assume that there exists an extension $\bar{h}\in \cU_\infty$ of $h$, \blue{with $\cU_\infty$ as in Proposition \ref{prop-assum}},  and replace $u$ by $\tilde{u}=u-\bar{h}$, $u_0$ by $\tilde{u}_0=u_0-\bar{h}(0)$, $g$ by $\tilde{g}:=g-D_t\bar{h}+\Delta \bar{h}$, $z$ by $\tilde{z}=z-\bar{h}$. After suppressing the tildes again, the model becomes
\begin{alignat}{3}
& D_tu - \Delta u + q(u+\bar{h}) + \phi(u+\bar{h})\psi(q) = g  \quad&&\mbox{ in }\Omega\times\I\label{ex1-1_a}\\
&u=0 &&\mbox{ on }\partial\Omega\times\I\\
& u(0) = u_0 &&\mbox{ in }\Omega\times\{0\}\label{ex1-2_a}
\end{alignat}
and the positivity condition on the exact state and its value under $\phi$ read as 
\[
z+\bar{h}\geq\underline{c}>0\,, \quad \phi(z+\bar{h})\geq0.
\]

Thus, we can use the spaces
\[U=H^1_0(\Omega),\quad Q=H^1(\Omega),\quad  H=L^2(\Omega)\]
and set 
\[
f(q,u)=-\Delta u + q(u+\bar{h}) + \phi(u+\bar{h})\psi(q).
\]

In the following, we verify Assumption \ref{ass-1}.
\begin{enumerate}
\item[\ref{A-extU}] 
\blue{\emph{Existence and uniqueness of exact state}
follows from Proposition \ref{prop-assum}.}

\item[\ref{A-coe}] \emph{Coercivity of $f$}
\begin{align*}
&\lan f(q,z)-f(q^*,z),q-q^*\ran_{Q^*,Q} \\
&\qquad= \lan (z+\bar{h})(q-q^*),q-q^*\ran_{Q^*,Q} + \lan \phi(z+\bar{h})(\psi(q)-\psi(q^*)),q-q^*\ran_{Q^*,Q}\\
&\qquad\geq \underline{c}\|q-q^*\|_H^2
=:\Ccoe\|q-q^*\|_H^2.
\end{align*}
Here, we invoke positivity of $z+\bar{h}$ and monotonicity of \textcolor{black}{$\phi(z+\bar{h})\psi$}.

\item[\ref{A-lip}] \emph{Local Lipschitz continuity of $f$}\\
We observe that with a constant $C$ depending on $M_z$, $C_\psi$, $\bar{h}$, and $\Omega$, one has
\begin{align*}
&\lan f(q,z)-f(q^*,z)- f'_q(\tilde{q},z)(q-q^*),v\ran_{U^*,U}\\
&=\lan \phi(z+\bar{h})(\psi(q)-\psi(q^*) - \psi'(\tilde{q})(q-q^*)),v\ran_{U^*,U}\\
&=\left\langle \phi(z+\bar{h})\left[\int_0^1 (\psi'(q^*+\lambda(q-q^*)) - \psi'(\tilde{q})) \,d\lambda(q-q^*)\right] ,v \right\rangle_{U^*,U}\\
&\leq CM_{z+\bar{h}}\|v\|_{L^6}\|q-q^*\|_{L^2}\left\|1+|q|^{\beta-1}+|q^*|^{\beta-1}+|\tilde{q}|^{\beta-1} \right\|_{L^3}\\
&\leq C_{H^1\to L^6}CM_{z+\bar{h}}\|v\|_{H^1}\|q-q^*\|_{L^2}\left(1+ \|q\|_{L^{3(\beta-1)}}^{\beta-1}+\|q^*\|_{L^{3(\beta-1)}}^{\beta-1}+\|\tilde{q}\|_{L^{3(2-\alpha)}}^{\beta-1} \right)\\
&=:L^{\tilde{q},z}(\|q\|_H)\|q-q^*\|_H\|v\|_U
\end{align*}
subjected to $1\leq \beta\leq 5/3$.\\
An example for this is $\alpha=3, \beta=5/3$ (see the remark below).
\item[\ref{A-C}] \emph{Choice of the linear operator $\cC$}\\
Taking $\cC(\|q\|_H)v= \left( L^{\tilde{q},z}(\|q\|_H)+ 1\right)(-\Delta) v,$ assumption \ref{A-C} holds.
\end{enumerate}

\begin{remark}
One of the feasible choices for the model nonlinearity is \textcolor{black}{$\phi(z)\psi(q^*)=z^3|q^*|^{\frac{2}{3}}q^*$}
with the exact coefficient $q^*$ \textcolor{black}{and exact state $z$} being nonnegative.
\end{remark}

\begin{remark}\label{remk-ex1-space}
The Lipschitz constant $L^{q_0,z}$ here depends on $M_z$, where $M_z$ depends on $T$ as shown in \eqref{ex0-boundz}. We thus have $L^{q_0,z}=L^{q_0,z,T}.$ This leads to $\Mtil=\Mtil^T, \Ntil=\Ntil^T$ in assumption \ref{A-C}, but still $C_{coe}=\underline{c}$ being independent of $T$ as $z\geq \underline{c}$ on $\I$. Therefore, comparing to Proposition \ref{prop-noisefree} we obtain for this example
\begin{alignat*}{3}
&u\in L^2(\I;U)\cap H^1(\T;U^*)\cap L^\infty(\I;H)&& \forall T>0,\\
&q \in L^2(\I;H)\cap H^1(\T;Q^*)\cap L^\infty(\I;H)\quad && \forall T>0.
\end{alignat*}
\end{remark}

\begin{remark}
As can be seen from the above estimate for \ref{A-coe}, coercivity could alternatively be achieved by imposing strict positivity of $\phi(z+\bar{h})$ (bounded away from zero) and uniform monotonicity of $\psi$.
\end{remark}

\begin{remark}
Smallness of $L^{\tilde{q},z}(\|q\|_H)$ could be obtained by assuming closeness of $q(t)$ and $\tilde{q}$ to $q^*$; as for $q(t)$, this can be bootstrapped from the exponential decay estimate \eqref{rate_exp}. This can be seen by rewriting $\psi'(q^*+\lambda(q-q^*)) - \psi'(\tilde{q})=\int_0^1\psi''(\tilde{q}+s(q^*+\lambda(q-q^*)-\tilde{q})) (q^*+\lambda(q-q^*)-\tilde{q})\, ds$ in the above estimate for \ref{A-lip}, which under the growth assumption $|\psi''(q)|\leq C_\psi(1+|q|^{\beta-2})$ gives an estimate that is qualitatively similar to the above -- including the conditions on $\alpha,\beta$ -- but with a constant that can be made small for $q(t)$ and $\tilde{q}$ close to $q^*$.
Preventing $L^{\tilde{q},z}(\|q\|_H)$, and therewith the factor in the definition of $\cC(\|q\|_H)$, from getting too large makes sense in order to avoid that the system \eqref{modelref-1}-\eqref{modelref-3} becomes very stiff, which might lead to high computational costs in simulation over a large time horizon.
\end{remark}

\subsection{Identification of a diffusion coefficient -- $a$ problem with nonlinear perturbation}
In this example, we study the estimation of the parameter $q$ in the model
\begin{align}
& D_tu - \nabla\cdot(q\nabla u) + \tilde{f}(u,q) = g  \quad\mbox{ in }\Omega\times\I\label{ex2-1}\\
& u(0) = u_0 \mbox{ in }\Omega\times\{0\}\label{ex2-2}
\end{align}
on a smooth bounded domain $\Omega\subseteq\mathbb{R}^3$, under the assumptions
\begin{align*}
&\tilde{f}(u,q):= u\Delta q + \phi(u)\psi(q)\\
&\text{the exact state } z(t,x)\leq -\underline{c}<0  \quad \forall x\in\Omega,\,t>0,\\
&\textcolor{black}{\phi(z)\psi(\cdot)} \text{ is monotone}, \quad\text{and}\quad
|\psi'(q)|\leq C_\psi(1+|q|^{\beta-1}).
\end{align*}
To achieve boundedness away from zero of $z$ by means of a maximum principle \cite[Lemma 2.1, Chapter 2]{Pao}, we impose inhomogeneous Dirichlet boundary conditions
\[
u=h\leq 0 \mbox{ on }\partial\Omega\times\I
\]
and assume that also 
\[
u_0\leq -\underline{c}\,, \quad g\leq -\underline{c}\Delta q^*\textcolor{black}{\,-M_zM_{q^*}} \,, \quad \Delta q^*\in L^\infty(\Omega)\,.
\]
This and the fact that $\hat{z}=z+\underline{c}$ solves
\begin{alignat*}{3}
& D_t\hat{z} - \nabla\cdot(q^*\nabla z) + \hat{z}\Delta q^* = g+\underline{c}\Delta q^*-\phi(z)\psi(q^*)\leq0  
\quad&&\mbox{ in }\Omega\times\I\\
&\hat{z}=h+\underline{c}\leq0 &&\mbox{ on }\partial\Omega\times\I\\
&\hat{z}(0) = u_0+\underline{c}\leq0 &&\mbox{ in }\Omega\times\{0\}
\end{alignat*}
enables the conclusion $\hat{z}\leq0$, yielding $z\leq -\underline{c}\leq 0$.
We again assume existence of an extension in order to work with homogeneous boundary conditions:
$\bar{h}\in \cU_\infty$ of $h$, \blue{with $\cU_\infty$ as in Proposition \ref{prop-assum}}, and replace $u$ by $\tilde{u}=u-\bar{h}$, $u_0$ by $\tilde{u}_0=u_0-\bar{h}(0)$, $g$ by $\tilde{g}:=g-D_t\bar{h}+\nabla\cdot(q\nabla\bar{h})-\bar{h}\Delta q$, $z$ by $\tilde{z}=z-\bar{h}$. After skipping the tildes, the model becomes
\begin{alignat}{3}
& D_tu - \nabla\cdot(q\nabla u) + u\Delta q + \phi(u+\bar{h})\psi(q)= g  \quad&&\mbox{ in }\Omega\times\I\\
&u=0 &&\mbox{ on }\partial\Omega\times\I\\
& u(0) = u_0 &&\mbox{ in }\Omega\times\{0\}.
\end{alignat}

According to the unique existence theory for linear parabolic PDEs, the diffusion coefficient must be positive; together with the continuous embedding $H^2(\Omega)\hookrightarrow L^\infty(\Omega)$, we choose the spaces
\[U=H_0^1(\Omega),\quad Q=H_0^1(\Omega)\cap H^2(\Omega),\quad  H=L^2(\Omega)\]
and denote the intermediate space $\Qtil:=H^1_0(\Omega)$; hence,  $Q\subset\Qtil\subseteq U$.\\ We then set 
\[
f(q,u)=- \nabla\cdot(q\nabla u) + u\Delta q + \phi(u+\bar{h})\psi(q).
\]
\textcolor{black}{$f$ is well-defined, since is it shown in \eqref{ex0-general} that for $q\in H^2(\Omega)$, thus $\Delta q\in L^2(\Omega)$, we have $f(q,u)\in U^*$.} \\

We now verify Assumption \ref{ass-1}.

\begin{enumerate}
\item[\ref{A-extU}] 
\blue{\emph{Existence and uniqueness of exact state} follows from Proposition \ref{prop-assum}.}

\item[\ref{A-coe}] \emph{Coercivity of $f$}
\begin{align} \label{EXa-coe}
&\lan f(q,z)-f(q^*,z),q-q^*\ran_{Q^*,Q} \nonumber\\
&= \lan - \nabla\cdot((q-q^*)\nabla z) + z\Delta (q-q^*),q-q^*\ran_{Q^*,Q} + \lan \phi(z+\bar{h})(\psi(q)-\psi(q^*)),q-q^*\ran_{Q^*,Q} \nonumber\\
&\geq \int_\Omega - \nabla\cdot((q-q^*)\nabla z)(q-q^*)\,dx + \int_\Omega z\Delta (q-q^*)(q-q^*)\,dx \\
&= \int_\Omega  (q-q^*)\nabla z\nabla(q-q^*)\,dx + \int_\Omega z\Delta (q-q^*)(q-q^*)\,dx \nonumber\\
&= -\int_\Omega  z\nabla\cdot((q-q^*)\nabla(q-q^*))\,dx + \int_\Omega z\Delta (q-q^*)(q-q^*)\,dx \nonumber\\
&= -\int_\Omega  z|\nabla(q-q^*)|^2\,dx \nonumber\\
&=:\Ccoe\|q-q^*\|_{\Qtil}^2. \nonumber
\end{align}
Above, we firstly invoke  monotonicity of \textcolor{black}{$\phi(z)\psi$}, then apply integration by parts while taking into account $q,q^*\in H^1_0(\Omega).$

As noted in Remark \ref{remk-strongnorm}, achieving the  $\Qtil$-norm here allows us to estimate the quantity $(q-q^*)$ with this strong norm in Assumption \ref{A-lip}. However, for the Lipschitz constant $L$, we need to stay with the weak norm $\|q\|_H$, since in \eqref{bound-et}-\eqref{bound-rt}, $L(\|q\|_{L^\infty(\I;H)})$ is required; this uniform boundedness in time is only attainable for $\|q(t)\|_H$, as proven in \eqref{rate_exp} or \eqref{bound-all}.

\item[\ref{A-lip}] \emph{Local Lipschitz continuity of $f$}\\
With a constant $C$ depending on $C_\phi$, $C_\psi,\Omega$, we see
\begin{align}\label{EXa-lip}
&\lan f(q,z)-f(q^*,z)- f'_q(\tilde{q},z)(q-q^*),v\ran_{U^*,U} \nonumber\\
&=\lan \phi(z+\bar{h})(\psi(q)-\psi(q^*) - \psi'(\tilde{q})(q-q^*)),v\ran_{U^*,U} \nonumber\\
&=\left\langle \phi(z+\bar{h})\left[\int_0^1 (\psi'(q^*+\lambda(q-q^*)) - \psi'(\tilde{q})) \,d\lambda(q-q^*)\right] ,v \right\rangle_{U^*,U} \nonumber\\
&\leq CM_{z+\bar{h}}\|v\|_{L^6}\|q-q^*\|_{L^6}\left\|1+|q|^{\beta-1}+|q^*|^{\beta-1}+|\tilde{q}|^{\beta-1} \right\|_{L^\frac{3}{2}} \nonumber\\
&\leq \left(C_{H^1\to L^6}\right)^2 CM_{z+\bar{h}}\|v\|_{H^1}\|q-q^*\|_{H_0^1}(1+\|z\|_{H^1}^\alpha)\left(1+ \|q\|_{L^\frac{3(\beta-1)}{2}}^{\beta-1}+\|q^*\|_{L^\frac{3(\beta-1)}{2}}+\|\tilde{q}\|_{L^\frac{3(\beta-1)}{2}} \right) \nonumber\\
&=:L^{\tilde{q},z}(\|q\|_H)\|q-q^*\|_{\Qtil}\|v\|_U,
\end{align}
provided that $1\leq\beta\leq 7/3$.\\
An example for this is $\alpha=3, \beta=7/3$ (see the remark below).
\item[\ref{A-C}] \emph{Choice of the linear operator $\cC$}\\
Taking $\cC(\|q\|_H)v= \left( L^{\tilde{q},z}(\|q\|_H)+ 1\right)(-\Delta) v,$ condition \ref{A-C} is fulfilled.
\end{enumerate}

\begin{remark}
One of the feasible choices for the model nonlinearity  is \textcolor{black}{$\phi(z)\psi(q^*)=z^3-z|q^*|^\frac{4}{3} q^*$}, with the exact coefficient $q^*$ nonnegative, its derivative $\Delta q^*$ bounded and the exact state $z$ nonpositive. The term $z^3$ is monotone w.r.t to $z$, and the term $-z|q^*|^\frac{4}{3} q^*$ is monotone w.r.t $q^*$ as $z<0$. In addition, $-z|q^*|^\frac{4}{3} q^*$ is linear in $z$ and hence plays the role of $cz$ with $c:=-|q^*|^\frac{4}{3} q^*\in L^2(\Omega)$ for $q^*\in Q$. In this setting, monotonicty in $z$ is not required.
\end{remark}

\begin{remark}
Similar to Remark \ref{remk-ex1-space} and together with the involvement of $\Qtil$, we obtain for this example
\begin{alignat*}{3}
&u\in L^2(\I;U)\cap H^1(\T;U^*)\cap L^\infty(\I;H)&& \forall T>0,\\
&q \in L^2(\I;\Qtil)\cap H^1(\T;\big(\Qtil\big)^*)\cap L^\infty(\I;H)\quad && \forall T>0.
\end{alignat*}
\end{remark}

\section{Conclusions and outlook}
In this paper, we have proposed and analyzed an online parameter identification method
for problems governed by nonlinear time-dependent PDEs. Our approach introduces a dynamic update law for both the state and the parameter via a model reference adaptive system. This system contains a linear PDE for the state and an auxiliary nonlinear one for the stationary parameter. Under suitable structural assumptions on the operators, unique existence results for this nonlinear adaptive system were given. In addition, by evaluating the error system, we proved that the solution to the adaptive system converges to the exact state and parameter for both exact as well as  noisy data. Our key contribution is that we tackle the case where the system of PDEs  is not only nonlinear with respect to the state, but also nonlinear with respect to the parameters; this situation has not been investigated thus far.

There are several future directions which could extend the study we have presented. Numerical experiments to illustrate the performance of our scheme are yet to be provided, and would be a natural continuation of this work. In the present work, we restrict ourselves to the case where the full state is measured. A modified scheme for the practically relevant situation of partial state observation is highly desirable. 
\blue{Likewise, perturbations not only of the data but also of the model $f$ would be of high practical interest.}
Finally, an extension of our result in the framework of time dependent parameters is also an appealing topic.

\section*{Acknowledgment} The work of BK was supported by the Austrian Science
Fund FWF under the grants P30054 and DOC 78.
Moreover, we wish to thank both reviewers for fruitful comments leading to an improved version of the manuscript.

\section*{Appendix}\label{appendix}
Before the problem \eqref{ex0-1}-\eqref{ex0-2} in Section \ref{sec:ex} is further discussed, we summarize the unique existence result for the dynamic equation
\begin{align}\label{InitProblem}
& D_tu(t) + f(q,u(t)) = g(t) \quad\text{for a.e.}\quad t\in(0,T], \qquad u(0) = u_0
\end{align}
in the same function space setting proposed in Section \ref{sec::adptive_sys}, except for the finite time interval.
The following auxiliary result is presented in the book by Roub\'i\v ceck \cite{Roubicek}:
\begin{auxiliary} \label{theo-Roubicek}
Assume that for fixed  $q\in Q$,
\begin{enumerate}[label=(S\arabic*)]
\item
$f$ is pseudomonotone
\item \label{S2}
$f$ is semi-coercive:
\[ \forall v \in U: \dupair{ f(q,v),v}_{U^*,U} \geq c_0|v|^2_U -c_1|v|_U-c_2\|v\|_H^2 \]
with some $c_0 > 0, \blue{c_1,c_2\in\R}$ and some seminorm $|.|_U$ satisfying $\|\cdot\|_U \leq c_{|\cdot|}(|\cdot|_U + \|\cdot\|_H)$ for some $c_{|.|}>0.$
\item \label{S3}
$f$ satisfies the regularity condition:
\begin{alignat*}{3}	
& g\in W^{1,\infty,2}(\T;U^*,U^*)\\
& u_0\in U \quad\text{such that}\quad f(u_0)-\blue{g(0)}\in H\\
& \lan f(q,u)-f(q,v), u-v\ran_{U^*,U}\geq  C_0|u-v|^2_U -C_2\|u-v\|_H^2
\end{alignat*}
with some $C_0>0, \blue{C_2\in\R}.$
\end{enumerate}
Then equation \eqref{InitProblem} has a unique solution $u \in W^{1,\infty,\infty}(\T;H,H)\cap W^{1,\infty,2}(\T;U,U),$ where $W^{1,p,q}(\T;U_1,U_2):=\{u\in L^p(\T;U_1):D_tu\in L^q(\T;U_2)\}$.
\proof \cite[Theorems 8.18, 8.31]{Roubicek}. \qed
\end{auxiliary}
\blue{Our plan is first verifying that \eqref{ex0-1}-\eqref{ex0-2} admits a solution $u$. In this step, the auxiliary result above will be used. Then, we lift the regularity of $u$ to the somewhat stronger space $L^\infty([0,T)\times\Omega)$. The second step  allows us to consider positivity of $u$, which is an important factor in the verification of Assumption \ref{ass-1}.}

\subsubsection*{\blue{Step 1: Unique existence of solution}}
Denoting in \eqref{ex0-1}\[-\nabla\cdot(a\nabla u) + cu + \phi(u)\psi(a,c):=f_1(u)+f_2(u)+f_3(u),\]
 it is evident that for the linear term,
\begin{enumerate}
\item 
$f_1:U\to U^*$ is monotone and continuous, thus pseudomonotone \cite[Lemma 6.7]{FRANCU} (one can also argue through coercivity).\\
$f_2:U\to U^*$ is strongly continuous, i.e., $u_n\rightharpoonup u$ implies $f_2(u_n)\to f_2(u)$. Indeed, for $u_n\rightharpoonup u$ in $U$, it holds that
\begin{align*}
&\|f_2(u_n)-f_2(u)\|_{U^*}=\sup_{\|v\|_U\leq 1}\int_\Omega c(u_n-u)v\,dx\\
&\leq \sup_{\|v\|_U\leq 1}\|c\|_{L^2(\Omega)}\|u_n-u\|^{\frac{1}{2}}_{L^6(\Omega)}\|u_n-u\|^{\frac{1}{2}}_{L^2(\Omega)}\|v\|_{L^6(\Omega)}\\
&\leq (C_{H^1\to L^6})^\frac{3}{2}\|c\|_{L^2(\Omega)}\|u_n-u\|^{\frac{1}{2}}_{U}\|u_n-u\|^{\frac{1}{2}}_{H} \to 0 \qquad \text{as }n\to\infty
\end{align*}
since the weakly convergent sequence $u_n$ is bounded in $U$, and the embedding $U=H^1(\Omega)\hookrightarrow L^2(\Omega)=H$ is compact. Thus, $f_2$ is pseudomonotone by \cite[Lemma 6.7]{FRANCU}.\\[2ex]
Hence, $f_{12}:=f_1+f_2$ is pseudomonotone \cite[Lemma 6.8]{FRANCU}.
\item $f_{12}$ is semi-coercive as required in \ref{S2}, since
\begin{align*}
&\langle f_{12}(u),u \rangle_{U^*,U}
=\int_\Omega (- \nabla\cdot(a\nabla u) + cu)u dx\\
&\geq \underline{a}\|\nabla u\|_{L^2(\Omega)}- C_{H^1\rightarrow L^6}\|c\|_{L^2(\Omega)}\left(\frac{\epsilon_1}{4\epsilon}+\epsilon\right)\|u\|_{H^1(\Omega)}^2 -  \frac{\blue{C_{H^1\rightarrow L^6}}}{16\epsilon\epsilon_1}\|c\|_{L^2(\Omega)}\|u\|^2_{L^2(\Omega)}, \\
&=:c_0 \|u\|_{U}^2 - c_2 \|u\|^2_{H}
\end{align*}
with $c_0>0$ as $\epsilon,\epsilon_1$ are arbitrarily small \cite[Section 3.1, (69)]{KNO}.
\item $f_{12}$ satisfies the regularity condition \ref{S3} with 
\[C_0:= c_0>0, \quad C_1=0, \quad C_2:= c_2\]
due to its linearity. We also assume, besides $g(0)\in H$, that
\begin{align}
u_0\in H^2(\Omega)
\end{align}
such that $f_{12}(u_0)-\blue{g(0)}\in H.$
\end{enumerate}

Based on this, unique existence of the solution to \eqref{ex0-1}-\eqref{ex0-2} then boils down the question of whether the nonlinear term $f_3$ fulfills the conditions in Auxiliary result \ref{theo-Roubicek}.\\
We make the assumption on monotonicity and growth of the linear term
\begin{align}
f_3=\phi(\cdot)\psi(a,c)\, \text{ is monotone,}\qquad\text{and}\qquad |\phi(u)|\leq C_\phi(1+|u|^\alpha)
\end{align}
then $f_3(u)\in U^*$, provided that
\begin{equation}\label{ex0-alpha1}
\begin{aligned}
&\alpha\leq 5 \qquad \text{if}\,\,\psi(a,c)\in L^\infty(\Omega),\quad \text{as}\,\, \|f_3(u)\|_{U^*}\leq C\|\psi(a,c)\|_{L^\infty(\Omega)}\|\phi(u)\|_{L^{6/5}},\\
&\alpha\leq 2 \qquad \text{if}\,\, \psi(a,c)\in L^2(\Omega),\quad\,\, \text{as}\,\, \|f_3(u)\|_{U^*}\leq C\|\psi(a,c)\|_{L^2(\Omega)}\|\phi(u)\|_{L^3}.
\end{aligned}
\end{equation}
We now observe
\begin{enumerate}
\item 
$f_{13}:=f_1+f_3$ is pseudomonotone referring to \cite[Arxiv version, Section 8.2, (S1)]{Nguyen}.\\[2ex]
Consequently, $f=f_{13}+f_2$ is pseudomonotone.
\item $f$ is semi-coercive, as
\begin{align*}
\lan f_3(u),u\ran_{U^*,U}\geq \lan f_3(0),u\ran_{U^*,U}\geq -C_\phi \|\psi(a,c)\|_H\|u\|_H \geq -C_{U\to H}C_\phi \|\psi(a,c)\|_H\|u\|_U,
\end{align*}
meaning $c_1:=C_{U\to H}C_\phi \|\psi(a,c)\|_H.$
\item
$f$ fulfills the regularity condition, since 
\[\lan f_3(u)-f_3(v),u-v\ran_{U^*,U}\geq 0.\]
As $u_0\in H^2(\Omega)\hookrightarrow L^\infty(\Omega)$, it follows that if $\psi(a,c)\in H,$ then \[\|f_3(u_0)-f_3(0)\|_H\leq C_\phi(2+(C_{H^2\to L^\infty})^\gamma\|u_0\|^\gamma_{L^\infty})\|\psi(a,c)\|_H.\] 
\end{enumerate}

This verification proves unique existence of the solution to the problem \eqref{ex0-1}-\eqref{ex0-2}
\begin{equation*}
u \in W^{1,\infty,\infty}(\T;H,H)\cap W^{1,\infty,2}(\T;U,U)\quad \forall T>0.
\end{equation*}

\subsubsection*{\blue{Step 2: Regularity of the solution}}
Now we turn to lift the regularity of $u$ to $L^\infty(\T;H^2(\Omega))\hookrightarrow L^\infty(\T\times\Omega)$ to attain its boundedness in time and space.\\
We observe from \eqref{ex0-1} that
\begin{alignat*}{3}
&D_tu\in L^\infty(0,T;L^2(\Omega)) \quad &&\text{as } u\in W^{1,\infty,\infty}(0,T;H,H) \\
& g\in L^\infty(0,T;L^2(\Omega)) &&\text{if we assume, together with } \eqref{ex0-general},\, g\in W^{1,\infty,2}(\T;H,U^*)\\
&\phi(u)\psi(a,c)\in L^\infty(0,T;L^2(\Omega)) \quad&&\text{if } \psi(a,c)\in L^\infty(\Omega), \text{and}\\
& &&\|\phi(u)\|_{L^2(\Omega)}\leq \|u\|_U^\gamma, \gamma\geq0 \text{ as } u\in W^{1,\infty,2}(\T;U,U),\\ 
& &&\text{thus constrain:}\quad \alpha \leq 3,
\end{alignat*}
which yields the remaining term  $-\nabla(a\nabla u)+cu=:b\in L^\infty(L^2(\Omega))$. If we assume $u=0$ on $\partial\Omega$, then $b$ also has zero boundary. Suppose that there exists $\tilde{c}\in L^\infty(\Omega)$ close to $c$, then for each $t\in(0,T)$,
\begin{align*}
&-\Delta u +\frac{\tilde{c}}{a}u = \frac{1}{a}(b+\nabla a\cdot\nabla u+(\tilde{c}-c)u )\\
&u- \left(-\Delta+\frac{\tilde{c}}{a}\right)^{-1} \left(\frac{1}{a}(\nabla a\cdot\nabla u+(\tilde{c}-c)u) \right)=: (\text{Id}-K)u=\left(-\Delta+\frac{\tilde{c}}{a}\right)^{-1}\left(\frac{1}{a}b\right),
\end{align*}
with $\left(-\Delta+\frac{\tilde{c}}{a}\right)^{-1}:L^2(\Omega)\to H^2(\Omega)$ for $\frac{\tilde{c}}{a}\in L^\infty(\Omega)$ and $\partial\Omega\in C^2$ being a bounded operator \cite[Section 6.3, Theorem 4]{Evans}.\\
Next, $K: H^2(\Omega)\cap H^1_0(\Omega)\to H^2(\Omega)\cap H^1_0(\Omega)$ is a linear bounded operator with $\|K\|<1$ if
\begin{align*}
\|Kv\|_{H^2(\Omega)}&=\left\|\left(-\Delta+\frac{\tilde{c}}{a}\right)^{-1}\left(\frac{1}{a}(\nabla a\cdot\nabla v +(\tilde{c}-c)v)\right)\right\|_{H^2(\Omega)}
\leq \frac{C^{\tilde{c},a}}{\underline{a}} \|\nabla a\cdot\nabla v +(\tilde{c}-c)v \|_{L^2(\Omega)}\\
&\leq \frac{C^{\tilde{c},a}}{\underline{a}} \left( C_{H^1\to L^6}\|\nabla a\|_{L^3(\Omega)}\|v \|_{H^2(\Omega)} + C_{H^2\to L^\infty}\|\tilde{c}-c\|_{L^2(\Omega)}\|v \|_{H^2(\Omega)} \right)\\
&<\|v \|_{H^2(\Omega)}
\end{align*}
for any $v\in H^2(\Omega)\cap H^1_0(\Omega)$, it holds 
\begin{align}
C_{H^1\to L^6}\|\nabla a\|_{L^3(\Omega)} + C_{H^2\to L^\infty}\|\tilde{c}-c\|_{L^2(\Omega)} < \underline{a}/C^{\tilde{c},a}.
\end{align}
Applying the Neumann series for $K$ with $\|K\|<1$, we have
\begin{align*}
\|u\|_{H^2(\Omega)}&=\left\|(\text{Id}-K)^{-1} \left(-\Delta+\frac{\tilde{c}}{a}\right)^{-1} \left(\frac{1}{a}b\right)\right\|_{H^2(\Omega)}\\
&\leq \frac{1}{1-\|K\|_{H^2(\Omega)\to H^2(\Omega)}}  \left\| \left(-\Delta+\frac{\tilde{c}}{a}\right)^{-1}\left(\frac{1}{a}b\right)\right\|_{H^2(\Omega)} \leq \frac{C^{\tilde{c},c}/\underline{a}}{1-\|K\|_{H^2(\Omega)\to H^2(\Omega)}} \| b\|_{L^2(\Omega)}.
\end{align*}
In conclusion, besides $u\in W^{1,\infty,\infty}(\T;H,H)\cap W^{1,\infty,2}(\T;U,U)$, we further have
\begin{equation*}
u \in  L^\infty(\T; H^2(\Omega))\hookrightarrow L^\infty(\T\times\Omega)\quad \forall T>0.
\end{equation*}
This result is summarized in Proposition \ref{prop-assum}.

\blue{
\begin{remark}
Alternatively, unique existence result for $u\in L^\infty([0,T)\times\Omega)$ can be constructed via a subsolutions supersolutions argument \cite[Section 1.5, Note and Comments 1.9]{Pao}. Monotonicity and growth of $u\mapsto\phi(u)\psi(a,c)$ remain as prerequisite assumptions, while conditions on $a, c, g$ might differ.
\end{remark}
}

\bibliographystyle{siam}
\nocite{*}
\bibliography{Reference}

\end{document}